\documentclass[12pt]{amsart}
\usepackage{graphicx}
\usepackage{subfigure}
\usepackage{latexsym}
\usepackage{amsmath}
\usepackage{amssymb}
\usepackage{mathrsfs}
\usepackage{amscd}
\usepackage{verbatim}

\pdfoutput=1

 \newtheorem{theorem}{Theorem}[section]
 
 \newtheorem{Prop}[theorem]{Proposition}
 \newtheorem{Lem}[theorem]{Lemma}
 \newtheorem{Cor}[theorem]{Corollary}
 
  \newtheorem{Coj}[theorem]{Conjecture}
  \newtheorem{Exa}[theorem]{Example}

\newcommand{\ba}{\begin{array}}
\newcommand{\ea}{\end{array}}
\newcommand{\beq}{\begin{equation}}
\newcommand{\eeq}{\end{equation}}

\newcommand\Let{\mathrel{\mathop:\!\!=}}

 \numberwithin{equation}{section}
 
\begin{document}

\title{ Connectedness of planar self-affine sets associated with non-consecutive collinear digit sets}

%\author{Ka-Sing Lau} \address{Department of Mathematics, The
%Chinese University of Hong Kong, Hong Kong}
%\email{kslau@@math.cuhk.edu.hk}

\author{King-Shun Leung} \address{ Department of Mathematics and Information Technology,
The Hong Kong Institute of Education, Hong Kong.}
\email{ksleung@@ied.edu.hk}

\author{Jun Jason Luo} \address{Department of Mathematics, The Chinese University of Hong Kong,
Hong Kong }
\email{jluo@@math.cuhk.edu.hk}

%\date{24 July 2010}
%
%
%\thanks {}
%
%
%\thanks{Supported by CNSF 19901025.}
%
%\subjclass{Primary: subject; Secondary: subject}
%\End Topmatter
%
%\date{\today}
%

\begin{abstract}
In the paper, we focus on the connectedness of planar self-affine
sets $T(A,{\mathcal{D}})$ generated by an integer expanding matrix
$A$ with $|\det (A)|=3$ and a collinear digit set
${\mathcal{D}}=\{0,1,b\}v$, where $b>1$ and $v\in {\mathbb{R}}^2$
such that $\{v, Av\}$ is linearly independent. We discuss the domain
of the digit $b$ to determine the connectedness of
$T(A,{\mathcal{D}})$. Especially, a complete characterization is
obtained when we restrict $b$ to be an integer. Some results on the
general case of $|\det (A)|> 3$ are obtained as well.
\end{abstract}

\maketitle

\bigskip

\begin{section} {\bf Introduction}
Let $M_n({\mathbb{Z}})$ denote the set of $n\times n$ matrices with
integer entries, let $A\in M_n({\mathbb{Z}})$ be expanding, i.e.,
all eigenvalues of $A$ have moduli strictly larger than $1$. Assume
$|\det (A)|=q$, and a finite set
${\mathcal{D}}=\{d_1,\dots,d_q\}\subset {\mathbb{R}}^n$ with
cardinality $q$, we call it a \emph{q-digit set}. It is well known
that there exists a unique \emph{self-affine set} $T\Let
T(A,{\mathcal{D}})$ \cite{LaWa} satisfying:
$$T= A^{-1}(T +
{\mathcal{D}})=\left\{\sum_{i=1}^{\infty}A^{-i}d_{j_i}: d_{j_i}\in
{\mathcal{D}}\right\}.$$ $T$ is called a \emph{self-affine tile} if
such $T$ has a nonvoid interior.

\medskip

The geometric and topological properties of $T(A,{\mathcal{D}})$
have been studied extensively. One of the interesting aspects is the
connectedness, in particular the disklikeness (i.e., homeomorphic to
a closed unit disc). It was asked by Gr\"ochenig and Haas
\cite{GrHa} that given an expanding integer matrix $A\in
M_n({\mathbb{Z}})$, whether there exists a digit set ${\mathcal{D}}$
such that $T(A,{\mathcal{D}})$ is a connected tile and they
partially solved the question in ${\mathbb R}^2$. Hacon et al.
\cite{HaSaVe} proved that any self-affine tile $T(A,{\mathcal{D}})$
with $2$-digit set is always pathwise connected. A systematical
study about this question was done by  \cite{KiLa}, \cite{KiLaRa},
\cite{HeKiLa} which mainly concerned  consecutive collinear digit
sets of the form $\{0,1,\dots,q-1\}v, v\in {\mathbb
Z}^n\setminus\{0\}$ via the algebraic property of the characteristic
polynomial of the matrix $A$. For some other case, Laarakker and
Curry \cite{LaCu} considered the connectedness of
$T(A,{\mathcal{D}})$ generated by an expanding matrix with rational
eigenvalues and a so-called centered canonical digit set.  In
${\mathbb R}^2$, the disklikeness is an interesting topic, Bandt and
Gelbrich \cite{BaGe}, Bandt and Wang \cite{BaWa}, and  Leung and Lau
\cite{LeLa} investigated the disklikeness of self-affine tiles in
terms of the neighbors of $T$. A translation of the tile $T+l$ is
called a neighbor of $T$ if $T\cap (T+l)\ne \emptyset$, where $l$ is
a lattice point (see Section 2 for details). Deng and Lau
\cite{DeLa}, and Kirat \cite{Ki} discussed the connectedness and
disklikeness of some other planar self-affine tiles with
non-collinear digit sets.

\medskip

Here we need to point out that, in \cite{KiLa} Kirat and Lau
obtained an interesting result in one dimensional space:

\medskip

\begin{Prop}\label{one dim. connectedness}
For $A=[q]$, and ${\mathcal{D}}=\{d_1,\dots,d_{|q|}\}\subset
{\mathbb{R}}$ with $q\in {\mathbb{Z}}, |q|\geq 2$, then
$T(A,{\mathcal{D}})$ is an interval (a connected tile) if and only
if, up to a translation, $ {\mathcal{D}}=\{0,a,\dots,(|q|-1)a\}$ for
some $a>0$.

\end{Prop}

\medskip

Unfortunately, this result cannot be extended to a higher
dimensional space in general. By using matrix expansion and
eigenvalue arguments,  Tan \cite{Ta} constructed a counterexample:

\medskip

\begin{Prop}\label{Tan's example}
Suppose $A\in M_2({\mathbb{Z}})$ with $|\det (A)|= 3$ is expanding,
whose characteristic polynomial is $f(x)=x^2-x-3$, and
${\mathcal{D}}=\{0,1,b\}v$, where $b>1$ and $v\in{\mathbb{R}}^2$
such that $\{v,Av\}$ is linearly independent. Then
$T(A,{\mathcal{D}})$ is connected if $8/5\leq b\leq 8/3$, and
$T(A,{\mathcal{D}})$ is disconnected if $b< \frac{\sqrt{13}-1}{2}$
or $b>\frac{\sqrt{13}+5}{2}$.
\end{Prop}

\medskip

As a generalization of Proposition \ref{one dim. connectedness}, we
can only have:

\medskip

\begin{theorem}\label{main result2}
Let $A\in M_n({\mathbb{Z}})$ be expanding with $|\det (A)|=q\geq 2$,
whose characteristic polynomial $f(x)=x^n\pm q $, the digit set is
${\mathcal{D}}=\{d_1,d_2,\dots,d_q\}v$ where $d_i\in {\mathbb{R}}$
and $v\in {\mathbb{R}}^n\setminus\{0\}$. Then the self-affine set
$T(A, {\mathcal D})$ is connected if and only if, up to a
translation, ${\mathcal{D}}$ is of the form
$\{0,1,2,\dots,(q-1)\}av$ for some $a>0$.
\end{theorem}

\medskip

However, more results on this kind of non-consecutive collinear
digit sets have not been obtained  yet.  Motivated by that, we try
to investigate this property in detail.

\medskip

For an expanding integer matrix $A\in M_2({\mathbb{Z}})$, it is
known by \cite{BaGe} that the characteristic polynomial of $A$ is
given by
$$f(x)=x^2+px+q, ~\text{with}~ |p|\leq q, ~\text{if}~ q\geq
2;\quad |p|\leq |q+2|, ~\text{if}~ q\leq -2.$$

When $|\det (A)| =3$, there are $10$ eligible characteristic
polynomials of $A$ as follows:
$$x^2\pm 3;\quad x^2\pm x+ 3;\quad x^2\pm 2x + 3;\quad x^2\pm 3x +
3;\quad x^2\pm x - 3.$$

Our main purpose in this paper is to study the self-affine sets
generated by $A$ with $|\det (A)|=3$ and ${\mathcal{D}}=\{0,1,b\}v$
with $b>1$ and $v\in{\mathbb{R}}^2$ such that $\{v,Av\}$ is linearly
independent. We get some criteria for $b$ to determine the
connectedness of $T(A,{\mathcal{D}})$. First if $b$ is an integer,
we have

\medskip

\begin{theorem}\label{main result3}
Let $|\det (A)|=3$, the characteristic polynomial of $A$ be
$f(x)=x^2+px\pm 3$ and ${\mathcal{D}}=\{0,1,m\}v$  where $2 \leq
m\in {\mathbb{Z}}$ such that $\{v,Av\}$ is linearly independent,.
Then we have:

(i) when $m=2$, $T(A, {\mathcal{D}})$ is always a connected tile;

(ii) when $m\geq 4$, $T(A, {\mathcal{D}})$ is always a disconnected
set;

(iii) when $m=3$, $T(A, {\mathcal{D}})$ is a connected set if
$f(x)=x^2\pm 2x+3 ~\text{or}~ x^2\pm 3x+3 ~\text{or}~ x^2\pm x-3$;
$T(A, {\mathcal{D}})$ is a disconnected set if $f(x)=x^2\pm 3
~\text{or}~ x^2\pm x+3$.

\end{theorem}

\medskip

Although there are many calculations in the proof, the method is
elementary. To get our desired results,  we apply the
Hamilton-Cayley theorem to find the potential neighbors of $T$ and
rule out the ineligible ones by a neighbor-generating algorithm. It
seems that this method can be applied to the more general situation.
By modifying the proof, we obtain

\medskip

\begin{theorem}\label{main results}
Assume that $A\in M_2({\mathbb{Z}})$ with $|\det (A)| =3$ is
expanding, a digit set ${\mathcal{D}}=\{0,1,b\}v$, where $b>1$ and $
v\in {\mathbb{R}}^2$ such that $\{v,Av\}$ is linearly independent.
Then we have:

Case 1: $f(x)= x^2\pm x + 3\quad T ~\text{is disconnected \quad if
\quad} b \geq 67/25 ~\text{or}~  b \leq 67/42;$

Case 2: $f(x)= x^2\pm 2x + 3\quad T ~\text{is disconnected \quad if
\quad}  b \geq 37/10 ~\text{or}~  b \leq 37/27;$

Case 3: $f(x)= x^2\pm 3x + 3\quad T ~\text{is disconnected \quad if
\quad} b \geq 33/10 ~\text{or}~  b \leq 33/23;$

Case 4: $f(x)= x^2\pm x - 3\quad T ~\text{is disconnected \quad if
\quad} b >19/5 ~\text{or}~ b<19/14.$
\end{theorem}

\medskip

It should be remarked that we mainly discuss each case under the
assumption of $b\geq 2$. If $1<b\leq 2$ then $ b/(b-1)\geq 2$,  with
the same argument we can consider the digit set
${\mathcal{D}}'=\{0,1,b/(b-1)\}(b-1)v$ which replaces the original
one,  since $T(A,{\mathcal{D}})$ is connected if and only if
$T(A,{\mathcal{D}}')$ is connected (just by a translation). So the
above related numbers come out.

\medskip

\begin{Coj}
We conjecture that there exists a critical value $c\geq 2$
(dependent on the characteristic polynomial of $A$) such that
$T(A,{\mathcal{D}})$ is connected if and only if $c/(c-1)< b\leq c$.
\end{Coj}

\medskip

Since our estimates about $b$ are rough, we can only solve the
conjecture partially. While if we suppose that $b$ is an integer,
then a positive answer is given as in Theorem \ref{main result3}.

\medskip

The rest of this paper is organized as follows. In Section 2, we
introduce some basic and useful results which play an important role
in our proofs. In Section 3, we first discuss the integer digit set
case where $b$ is an integer and then prove Theorems \ref{main
result2} and \ref{main result3}. Theorem \ref{main results} will be
shown in Section 4. In the last section we consider the general
planar self-affine sets with $|\det (A)|>3$ and non-consecutive
collinear digit sets, some sufficient conditions for $T(A,{\mathcal
D})$ to be connected are obtained.

\end{section}

\bigskip

\begin{section}{\bf Preliminaries}

In this section, we prepare some elementary knowledge about the
geometric properties of self-affine sets which will be used
frequently in the paper.

\medskip

We define
$${\mathcal{E}}=\{(d_i,d_j):(T+d_i)\cap(T+d_j)\ne\emptyset,
d_i,d_j\in {\mathcal{D}}\}$$ to be the set of edges for the digit
set ${\mathcal{D}}$. We say that $d_i$ and $d_j$ are
$\mathcal{E}$-connected if there exists a finite sequence (or path)
$\{d_{j_1},\dots,d_{j_k}\}\subset {\mathcal{D}}$ such that
$d_i=d_{j_1},d_j=d_{j_k}$ and $(d_{j_l},d_{j_{l+1}})\in
{\mathcal{E}}, 1\leq l \leq k-1.$

\medskip

It is easy to check that $(d_i,d_j)\in {\mathcal{E}}$ if and only if
$d_i-d_j=\sum_{k=1}^{\infty}A^{-k}v_k$ where $v_k\in
\Delta\mathcal{D}\Let \mathcal{D}-\mathcal{D}$. Then we get a
criterion of connectedness of a self-affine set by using a graph
argument on $\mathcal{D}$:

\medskip

\begin{Prop}(\cite{KiLa}) \label{e-connected prop}
A self-affine set $T$ with a digit set $\mathcal{D}$ is connected if
and only if any two $d_i, d_j\in {\mathcal{D}}$  are
$\mathcal{E}$-connected.
\end{Prop}

\medskip

When the digit set ${\mathcal{D}}$ can be written as
$\{d_1v,\dots,d_qv\}$ for some non-zero vector $v\in {\mathbb{R}}^n$
and $d_1< d_2<\cdots < d_q, ~d_i\in {\mathbb{R}}$, ${\mathcal{D}}$
is said to be collinear. If $d_{i+1}-d_i= constant (=1)$,
${\mathcal{D}}$ is called \emph{consecutive collinear digit set}.
Let $D=\{d_1,\dots,d_q\}, \Delta D= D-D=\{d=d_i-d_j: d_i, d_j\in
D\}$. Then ${\mathcal{D}}=Dv$ and ${\Delta\mathcal{D}}=\Delta Dv$.
It is easy to see that the connectedness of $T$ is invariant under a
translation of the digit set, hence we always assume that $d_1=0$.
The radix expansion of a point $x=\sum_{i=1}^{\infty}a_iA^{-i}v \in
{\mathbb R}^2$ is given by $0.a_1a_2a_3\ldots$ An overbar denotes
repeating digits as in $0.12\overline{301}=0.12301301301\ldots$
Likewise, $a_{-2}a_{-1}a_0.a_1a_2a_3\ldots$ represents a point
$$a_{-2}A^2v+a_{-1}Av+a_0v+\sum_{i=1}^{\infty}a_iA^{-i}v.$$

Let ${\mathbb{Z}}[x]$ denote the class of polynomials with integer
coefficients. We call a polynomial $f(x)\in {\mathbb{Z}}[x]$ an {\it
expanding polynomial} if all its roots have moduli strictly bigger
than $1$. Note that a matrix $A\in M_n({\mathbb{Z}})$ is expanding
if and only if its characteristic polynomial is expanding.

\medskip

We say that a monic polynomial $f(x)\in {\mathbb{Z}}[x]$ with
$|f(0)|=q$ has the \emph{Height Reducing Property} (HRP) if there
exists $g(x)\in {\mathbb{Z}}[x]$ such that
$$g(x)f(x)=x^k+a_{k-1}x^{k-1}+\cdots+a_1x\pm q,$$ where $|a_i|\leq
q-1, i=1,\dots,k-1.$

\medskip

This property was introduced by Kirat and Lau \cite{KiLa} to study
the connectedness of self-affine tiles associated with consecutive
collinear digit sets. It was proved that:

\medskip

\begin{Prop} \label{HRP}
Let $A\in M_n({\mathbb{Z}})$ with $|\det (A)|=q$ be expanding and
let ${\mathcal{D}}=\{0,1,2,\dots,(q-1)\}v$  be a consecutive
collinear digit set in ${\mathbb{R}}^n$.  Suppose the characteristic
polynomial $f(x)$ of $A$ has HRP, then $T$ is connected.
\end{Prop}

\medskip

In \cite{KiLaRa}, Kirat et al conjectured that all expanding integer
monic polynomials have HRP.  Akiyama and Gjini \cite{AkGj} solved it
up to degree $4$. However it is still open for arbitrary degree.
Recently, He et al \cite{HeKiLa} developed an algorithm to check HRP
for any monic polynomial, and Akiyama et al \cite{AkDrJa} studied
HRP of algebraic integers on canonical number systems from another
perspective.

\medskip

In the sequel, we always consider the planar self-affine set
$T(A,{\mathcal{D}})$ associated with a collinear digit set of the
form $\{d_1v,\dots,d_qv\}$, where $v$ is a vector in
${\mathbb{R}}^2$ such that $\{v, Av\}$ is linearly independent. It
is known that $\{v,Av\}$ is always linearly independent provided
$|\det (A)|=3$ (or any prime number) and $v\in {\mathbb
Z}^2\setminus \{0\}$. Denote the characteristic polynomial of $A$ by
$f(x)=x^2+px+q$, where $p,q \in {\mathbb Z}$. Then we can regard $A$
as the companion matrix of $f(x)$ for simplicity, i.e.,
$$A=\left[
        \begin{array}{rrrr}
          0  & -q \\
           1 &  -p   \\
    \end{array}
    \right].$$
Let $\Delta=p^2-4q$ be the discriminant. Define $\alpha_i,\beta_i$
by $$A^{-i}v=\alpha_iv+\beta_iAv, \quad i=1,2,\dots$$

Applying the Hamilton-Cayley theorem $f(A)=A^2+pA+qI=0$, where $I$
is a $2\times 2$ identity matrix, and the definitions of
$\alpha_i,\beta_i$ above, we have a lemma:

\medskip

\begin{Lem}(\cite{Le}) \label{evaluation}
Let $\alpha_i,\beta_i$ be defined as the above. Then
$q\alpha_{i+2}+p\alpha_{i+1}+\alpha_i=0$ and
$q\beta_{i+2}+p\beta_{i+1}+\beta_i=0$, i.e.,
$$\left[
        \begin{array}{rr}
          \alpha_{i+1} \\
          \alpha_{i+2} \\
    \end{array}
    \right]=\left[
        \begin{array}{rrrr}
          0  & 1 \\
           -1/q &  -p/q   \\
    \end{array}
    \right]^i \left[
        \begin{array}{rr}
          \alpha_1 \\
          \alpha_2 \\
    \end{array}
    \right]; \quad \left[
        \begin{array}{rr}
          \beta_{i+1} \\
          \beta_{i+2} \\
    \end{array}
    \right]=\left[
        \begin{array}{rrrr}
          0  & 1 \\
           -1/q &  -p/q   \\
    \end{array}
    \right]^i \left[
        \begin{array}{rr}
          \beta_1 \\
          \beta_2 \\
    \end{array}
    \right]
$$
and $\alpha_1=-p/q, \alpha_2=(p^2-q)/q^2; \beta_1=-1/q,
\beta_2=p/q^2$. Moreover for $\Delta\ne 0$, we have
$$\alpha_i=\frac{q(y_1^{i+1}-y_2^{i+1})}{\Delta^{1/2}}
\quad\text{and}\quad \beta_i=\frac{-(y_1^i-y_2^i)}{\Delta^{1/2}},$$
where $y_1=\frac{-p+\Delta^{1/2}}{2q}$ and
$y_2=\frac{-p-\Delta^{1/2}}{2q}$ are the two roots of $qx^2+px+1=0$.
\end{Lem}

\medskip

Since $f(x)=x^2+px+q$ is expanding, the two roots $y_1, y_2$ of
$x^2f(x^{-1})=qx^2+px+1$ have moduli less than $1$. It follows that
the following two series converge:
\begin{equation*}
\tilde{\alpha}\Let
\sum_{i=1}^{\infty}|\alpha_i|,\quad
\tilde{\beta}\Let
\sum_{i=1}^{\infty}|\beta_i|.
\end{equation*}

From the previous lemma, it is easy to say, when $\Delta< 0$ (hence
$q\geq 2$), then
$$|\alpha_i|\leq \frac{2q|y_1^{i+1}|}{|\Delta^{1/2}|}= \frac{2q^{-(i-1)/2}}{(4q-p^2)^{1/2}}
\quad\text{and}\quad|\beta_i|\leq\frac{2|y_1^i|}{|\Delta^{1/2}|}=\frac{2q^{-i/2}}{(4q-p^2)^{1/2}}.$$

Hence, the upper bounds of $\tilde{\alpha}, \tilde{\beta}$ are
estimated by:
\begin{eqnarray}
\tilde{\alpha}  & \leq & \sum_{i=1}^{n-1}|\alpha_i|+
\frac{2q^{-(n-1)/2}}{(1-q^{-1/2})(4q-p^2)^{1/2}}, \label{estimates of two series1}\\
\tilde{\beta}  &\leq & \sum_{i=1}^{n-1}|\beta_i|+
\frac{2q^{-n/2}}{(1-q^{-1/2})(4q-p^2)^{1/2}}. \label{estimates of
two series2}
\end{eqnarray}
Since $q\geq 2$, we can find very accurate upper estimates of
$\tilde{\alpha}$ and $ \tilde{\beta}$ by taking $n=13$ or any larger
integer. This is the most important tool in our proofs.

\medskip

Let $L:=\{\gamma v+\delta Av: \gamma,\delta\in {\mathbb{Z}}\}$ be
the \emph{lattice} generated by $\{v,Av\}$. For $l\in
L\setminus\{0\},~ T+l$ is called a \emph{neighbor} of $T$ if
$T\cap(T+l)\ne \emptyset.$ It is easy to show that $T+l$ is a
neighbor of $T$ if and only if $l$ can be expressed in the form:
$$l=\sum_{i=1}^{\infty}b_iA^{-i}v \in T-T, ~\text{where}~ b_i\in \Delta D.$$

After a few calculations, we have

\medskip

\begin{Lem}(\cite{Le}) \label{neighbor generator}
If $T+l$ is a neighbor of $T$, where $l=\gamma v+\delta
Av=\sum_{i=1}^{\infty}b_iA^{-i}v$, then   $|\gamma|\leq
\max_i|b_i|\tilde \alpha$ and  $|\delta|\leq \max_i|b_i|\tilde
\beta$.  Moreover, this implies that $T+l^{\prime}$ is another
neighbor of $T$ satisfying
$l^{\prime}=Al-b_1v=\gamma^{\prime}v+\delta^{\prime}Av$ with
$\gamma^{\prime}=-(q\delta+b_1)$ and
$\delta^{\prime}=\gamma-p\delta$.
\end{Lem}

\medskip

By repeatedly using Lemma \ref{neighbor generator}, we can construct
a sequence of neighbors of $T$: $\{T+l_n\}_{n=0}^{\infty}$, where
$l_0=l$ and $l_n=\gamma_n v+\delta_n Av, n\geq 1$ by the following
inductive formula:
\begin{equation}\label{inductive formula of gamma_n and delta_n}
\left[
        \begin{array}{rr}
          \gamma_n \\
          \delta_n \\
    \end{array}
    \right]=A^n \left[
        \begin{array}{rr}
          \gamma \\
          \delta \\
    \end{array}
    \right]-\sum_{i=1}^n A^{i-1}\left[
        \begin{array}{rr}
          b_{n+1-i} \\
          0\\
    \end{array}
    \right].
\end{equation}
Moreover, $|\gamma_n|\leq \max_i|b_i|\tilde \alpha$ and
$|\delta_n|\leq \max_i|b_i|\tilde \beta$ hold for any $n\geq 0$.

\medskip

\begin{Lem}
Let $T_1=T(A,{\mathcal{D}})$ and $T_2=T(-A,{\mathcal{D}})$. Then
$T_1+l$ is a neighbor of $T_1$ if and only if $T_2+l$ is a neighbor
of $T_2$.
\end{Lem}

\medskip

\begin{proof}
If $l\in T_1-T_1$, then
$$l=\sum_{i=1}^{\infty}b_iA^{-i}v=\sum_{i=1}^{\infty}b_{2i}(-A)^{-2i}v+\sum_{i=1}^{\infty}(-b_{2i-1})(-A)^{-2i+1}v,$$ i.e., $l \in
T_2-T_2 $  and vice versa.
\end{proof}

\medskip

An immediate corollary of the lemma follows:

\medskip

\begin{Cor}\label{symmetry of char. poly.}
If the characteristic polynomial of the expanding matrix $A$ is
$x^2+px+q$ and that of $B$ is $x^2-px+q$. Then the self affine set
$T(A,{\mathcal{D}})$ is connected if and only if
$T(B,{\mathcal{D}})$ is connected.
\end{Cor}

\end{section}

\bigskip

\begin{section} {\bf Integer collinear digit set: ${\mathcal{D}}=\{0,1,m\}v, 2 \leq m\in {\mathbb Z}.$}

In the section, we mainly consider the connectedness of the
self-affine set $T(A,{\mathcal{D}})$ associated with a special digit
set ${\mathcal{D}}=\{0,1,m\}v$, where $2\leq m\in {\mathbb{Z}}$.

\medskip

Recall that, for $|\det (A)| =3$, there are $10$ eligible
characteristic polynomials of $A$ of the form: $f(x)=x^2+px\pm 3$.
We are going to discuss them separately.

\medskip

\subsection{ When $p=0$} We first provide a general result:

\begin{theorem}\label{thm of p=0}
Let $A\in M_n({\mathbb{Z}})$ be expanding with $|\det (A)|=q\geq 2$,
whose characteristic polynomial $f(x)=x^n\pm q $, the digit set is
${\mathcal{D}}=\{d_1,d_2,\dots,d_q\}v$ where $d_i\in {\mathbb{R}}$
and $v\in {\mathbb{R}}^n\setminus\{0\}$. Then the self-affine set
$T(A, {\mathcal D})$ is connected if and only if, up to a
translation, ${\mathcal{D}}$ is of the form
$\{0,1,2,\dots,(q-1)\}av$ for some $a>0$.
\end{theorem}

\medskip

\begin{proof}
Consider $f(x)=x^n-q$ only.  Since $f(A)=0$, we have
$A^{-n}=q^{-1}I$. Let $y=\sum_{i=1}^{\infty}a_iA^{-i}v\in T$, where
$a_i\in D=\{d_1,d_2,\dots,d_q\}$. Then
$$y=\sum_{k=0}^{n-1}(\sum_{j=1}^{\infty}a_{jn-k}A^{-jn})A^kv=\sum_{k=0}^{n-1}(\sum_{j=1}^{\infty}a_{jn-k}q^{-j})A^kv.$$

If $T$ is connected, then for every $k\in\{1,2,\dots, n-1\}$, along
the direction $A^kv$, the coordinate set
$$\{\sum_{j=1}^{\infty}a_{jn-k}q^{-j}: a_{jn-k}\in D\}$$ is an
interval, which is equivalent to say $D$ is a translation of
$\{0,1,2,\dots,(q-1)\}a$ for some $a>0$ by Proposition \ref{one dim.
connectedness}. Consequently, ${\mathcal{D}}$ is of the form
$\{0,1,2,\dots,(q-1)\}av$.

The converse is also true by Proposition \ref{HRP}, since $f(x)$
above has HRP.
\end{proof}

\medskip

\begin{Cor}
Assume that $A\in M_2({\mathbb{Z}})$ with $|\det (A)| =3$ is
expanding, its characteristic polynomial is $f(x)=x^2\pm 3$,
${\mathcal{D}}=\{0,1,b\}v$ is a collinear digit set, where $1< b\in
{\mathbb R}$ and $v\in {\mathbb{R}}^2\setminus\{0\}$. Then
$T(A,{\mathcal{D}})$ is connected if and only if $b=2$.
\end{Cor}

\medskip
\subsection{ When $p\ne 0$ and $m\geq 4$}

We will show that, in all cases, $(T+v)\cap (T+mv)=\emptyset$ and
$T\cap (T+mv)=\emptyset$. Hence by Proposition \ref{e-connected
prop}, $T$ is disconnected. (see Figure \ref{figure 1} where we take
the vector $v=\left[\begin{array}{rr} 1\\0\end{array}\right]$)

\medskip

\begin{figure}[h]
  \centering
  \subfigure[ $x^2+x+3$]{
  \includegraphics[width=5cm]{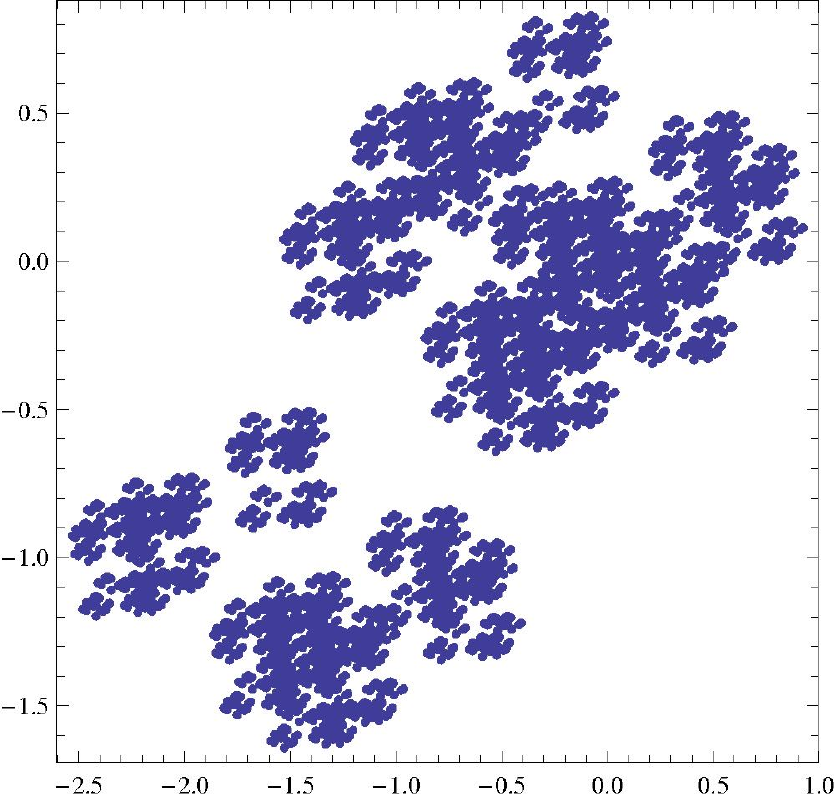}
 }
\qquad
  \subfigure[$x^2+2x+3$]{
   \includegraphics[width=5cm] {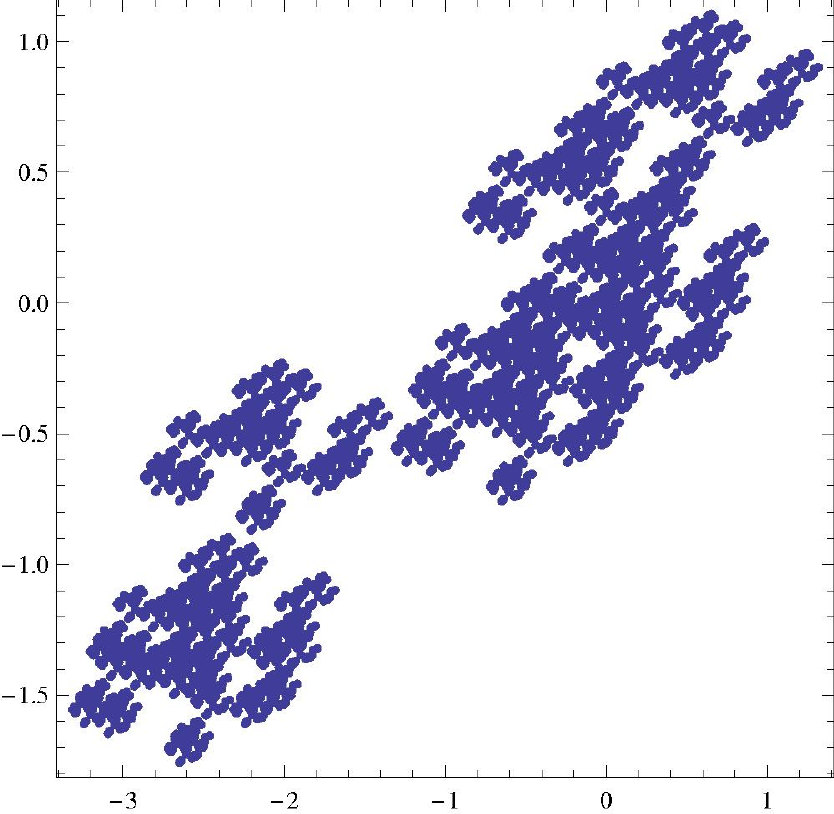}
 }\\
 \subfigure[$x^2+3x+3$]{
   \includegraphics[width=5cm]{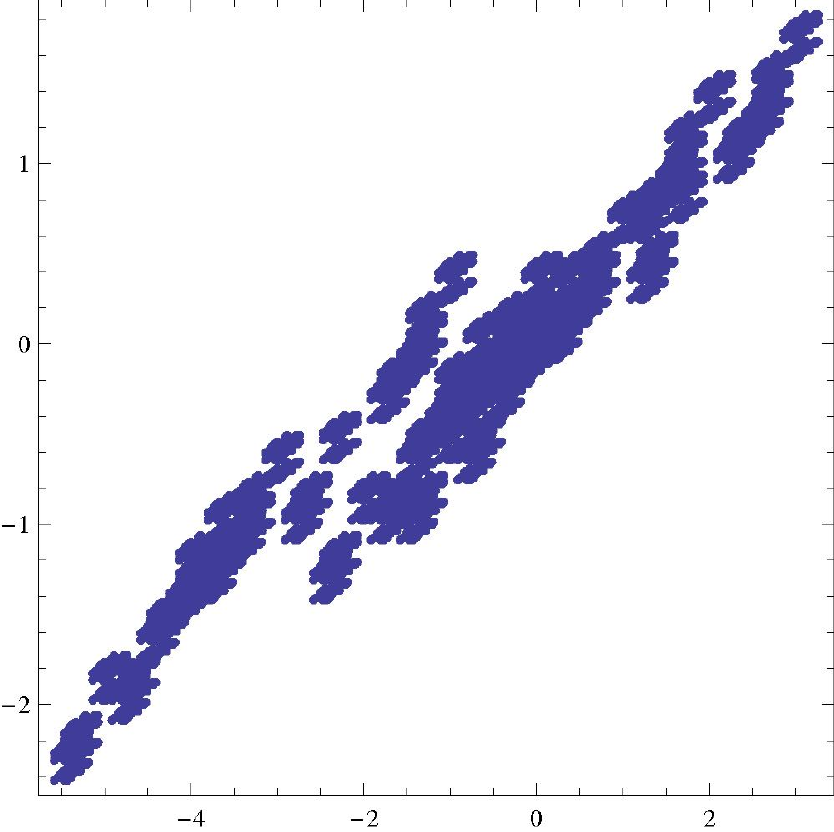}
 }
\qquad
 \subfigure[$x^2+x-3$]{
   \includegraphics[width=5cm]{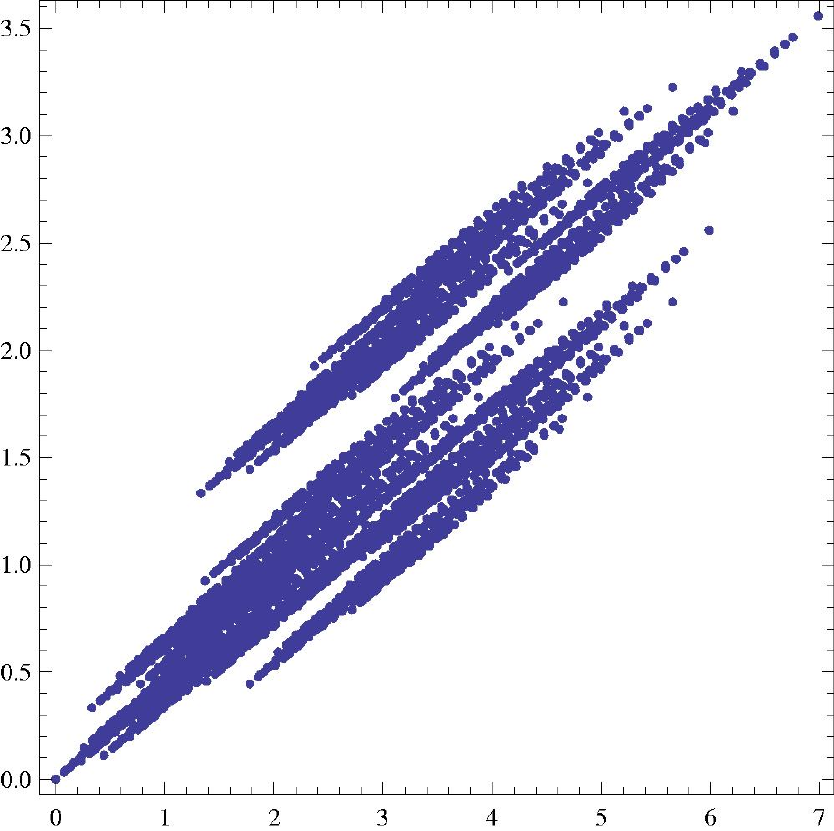}
 }
\caption{ When $p\ne 0, m=4$.} \label{figure 1}
\end{figure}

\begin{theorem}\label{main thm. of special digit sets}
Let the characteristic polynomial of $A$ be $f(x)=x^2+px\pm 3 ~(p\ne
0)$ and ${\mathcal{D}}=\{0,1,m\}v, m\geq 4$ be the digit set. Then
$T(A, {\mathcal{D}})$ is disconnected.
\end{theorem}

\medskip

\begin{proof}
It suffices to prove the cases for $p>0$ by Corollary \ref{symmetry
of char. poly.}.

\medskip

{\bf Case $1: f(x)=x^2+x+3$.} Let $T+l$ be a neighbor of $T$, then
$l=\gamma v+\delta Av=\sum_{i=1}^{\infty}b_iA^{-i}v$, $b_i\in \Delta
D=\{0,\pm 1,\pm(m-1),\pm m\}$.  By Lemma \ref{evaluation} and
(\ref{estimates of two series1}), (\ref{estimates of two series2}),
we have
\begin{align}\label{estimate of gamma and delta 1}
\tilde{\alpha}<0.88;\quad \tilde{\beta}<0.63.
\end{align}
Then $|\gamma|\leq m\tilde{\alpha}<0.88m,~ |\delta|\leq
m\tilde{\beta}<0.63m.$

Suppose $(T+v)\cap (T+mv)\ne\emptyset$, then
$(m-1)v=\sum_{i=1}^{\infty}b_iA^{-i}v, b_i\in \Delta D$. By
(\ref{inductive formula of gamma_n and delta_n}),  we obtain  $l_1=
-(3(m-1)+b_2)v-(m-1+b_1)Av$,  where $|3(m-1)+b_2|<0.88m$. Since
$|3(m-1)+b_2|\geq 3(m-1)-m=2m-3$, it follows $m<3/1.12$, which
contradicts $m\geq 4$.

Suppose $T\cap (T+mv)\ne\emptyset$, similarly, we get $l_1=
-(3m+b_2)v-(m+b_1)Av$, where $|3m+b_2|<0.88m$. Since $|3m+b_2|\geq
3m-m=2m$, a contradiction follows.

\medskip

{\bf Case $2: f(x)=x^2+2x+3$.} Let $T+l$ be a neighbor of $T$, then
$l=\gamma v+\delta Av=\sum_{i=1}^{\infty}b_iA^{-i}v$, where $b_i\in
\Delta D$. By Lemma \ref{evaluation} and (\ref{estimates of two
series1}), (\ref{estimates of two series2}), we have
\begin{align}\label{estimate of gamma and delta 2}
\tilde{\alpha}<1.17, \quad \tilde{\beta}<0.73.
\end{align}
Then $|\gamma|\leq m\tilde{\alpha}< 1.17m, ~ |\delta|\leq
m\tilde{\beta}< 0.73m$.

Suppose $(T+v)\cap (T+mv)\ne\emptyset$, then
$(m-1)v=\sum_{i=1}^{\infty}b_iA^{-i}v, b_i\in \Delta D$. By
(\ref{inductive formula of gamma_n and delta_n}),  we obtain
$l_1=-(3(m-1)+b_2)v-(2(m-1)+b_1)Av$. By (\ref{estimate of gamma and
delta 2}), then $2m-3=3m-3-m\leq |3(m-1)+b_2|< 1.17m$, it follows
that $m< 4$ which contradicts the assumption $m\geq 4$.

Suppose $T\cap (T+mv)\ne\emptyset$, similarly,  by (\ref{inductive
formula of gamma_n and delta_n}), we find a neighbor
$T-(3m+b_2)v-(2m+b_1)Av$. So $2m=3m-m\leq 3m+b_2\leq |3m+b_2|<
1.17m$, then $m<0$, a contradiction follows.

\medskip

{\bf Case $3: f(x)=x^2+3x+3$.} Let $T+l$ be a neighbor of $T$, then
$l=\gamma v+\delta Av=\sum_{i=1}^{\infty}b_iA^{-i}v$, where $b_i\in
\Delta D$. By Lemma \ref{evaluation} and (\ref{estimates of two
series1}), (\ref{estimates of two series2}), we have
\begin{align}\label{estimate of gamma and delta 3}
\tilde{\alpha}<2.24, \quad \tilde{\beta}<1.08.
\end{align}
Then $|\gamma|\leq m\tilde{\alpha}< 2.24m, \quad |\delta|\leq
m\tilde{\beta}< 1.08m$.

Suppose $T\cap (T+mv)\ne\emptyset$, then
$mv=\sum_{i=1}^{\infty}b_iA^{-i}v, b_i\in \Delta D$. By
(\ref{inductive formula of gamma_n and delta_n}), we obtain
$l_1=-(3m+b_2)v-(3m+b_1)Av$, where $|3m+b_1|< 1.08m$. Then $b_1<
-1.92m$, which is not in $\Delta D$.  This is ridiculous.

Suppose $(T+v)\cap (T+mv)\ne\emptyset$, similarly, by
(\ref{inductive formula of gamma_n and delta_n}),
$l_2=(9m-9+3b_1-b_3)v+(6m-6+3b_1-b_2)Av$, where
$|9m-9+3b_1-b_3|<2.24m$. Since $|9m-9+3b_1-b_3|\geq 5m-9$, it
follows $m<9/2.76<4$, which contradicts the assumption $m\geq 4$.

\medskip

{\bf Case $4: f(x)=x^2+x-3$.} In this case, the discriminant
$\Delta>0$. By Lemma \ref{evaluation}, we have $\alpha_1=1/3,
~\alpha_2=4/9, ~\beta_1=1/3, ~\beta_2=1/9$, and
$$\left[
        \begin{array}{rr}
          \alpha_{i+1} \\
          \alpha_{i+2} \\
    \end{array}
    \right]=B^i \left[
        \begin{array}{rr}
          1/3 \\
          4/9 \\
    \end{array}
    \right]; \quad \left[
        \begin{array}{rr}
          \beta_{i+1} \\
          \beta_{i+2} \\
    \end{array}
    \right]=B^i \left[
        \begin{array}{rr}
          1/3\\
          1/9 \\
    \end{array}
    \right],
$$ where
$B=\left[\begin{array}{cc}
0 & 1\\
1/3 & 1/3
\end{array}\right]$.
Hence $\alpha_i\geq 0$ and $\beta_i\geq 0$ for all $i$. It yields
that $$\tilde{\alpha}=\left[
        \begin{array}{cc}
          1 & 0
          \end{array}
    \right]\sum_{i=0}^{\infty}B^i \left[
        \begin{array}{rr}
          1/3 \\
          4/9 \\
    \end{array}
    \right]=\left[
        \begin{array}{cc}
          1 & 0
          \end{array}
    \right](I-B)^{-1}\left[
        \begin{array}{rr}
          1/3 \\
          4/9 \\
    \end{array}
    \right]=2.$$

Similarly, we also have $\tilde\beta=1$.

Let $T+l$ be a neighbor of $T$, then $l=\gamma v+\delta
Av=\sum_{i=1}^{\infty}b_iA^{-i}v$, where $b_i\in \Delta D$. By Lemma
\ref{neighbor generator},
\begin{align}\label{estimate of gamma and delta 4}
|\gamma|\leq m\tilde{\alpha}=2m, \quad |\delta|\leq
m\tilde{\beta}=m.
\end{align}

Suppose $(T+v)\cap (T+mv)\ne\emptyset$, then
$(m-1)v=\sum_{i=1}^{\infty}b_iA^{-i}v, b_i\in \Delta D$. By
(\ref{inductive formula of gamma_n and delta_n}), $l_1=
-(3m-3+3b_1+b_3)v+(4m-4+b_1-b_2)Av$,  where $|4m-4+b_1-b_2|\leq m$.
Since $|4m-4+b_1-b_2|\geq 2m-4$, it follows that $m=4,b_1=-4,b_2=4$.
Hence  $l_1=-(b_3-3)v+4Av$. Using (\ref{inductive formula of gamma_n
and delta_n}) again,  we get $l_2=(12-b_4)v-(1+b_3)Av$. Then $b_4=4$
from $|12-b_4|\leq 8$, and $l_2=8v-(1+b_3)Av$. Repeatedly using
(\ref{inductive formula of gamma_n and delta_n}), then
$l_3=-(3+3b_3+b_5)v+(9+b_3)Av$, where $|9+b_3|\leq 4$.  Since
$|9+b_3|\geq 5$, a contradiction follows.

Suppose $T\cap (T+mv)\ne\emptyset$, similarly,
$T-(3m+3b_1+b_3)v+(4m+b_1-b_2)Av$ is a neighbor of $T$ by
(\ref{inductive formula of gamma_n and delta_n}), where
$|4m+b_1-b_2|\leq m$. Since $|4m+b_1-b_2|\geq 2m$, it is impossible.

Therefore, in all cases, $T(A,{\mathcal{D}})$ is disconnected and
the theorem follows.

\end{proof}

\medskip

\begin{figure}[h]
 \centering
  \subfigure[ $x^2+x+3$]{
 \includegraphics[width=5cm]{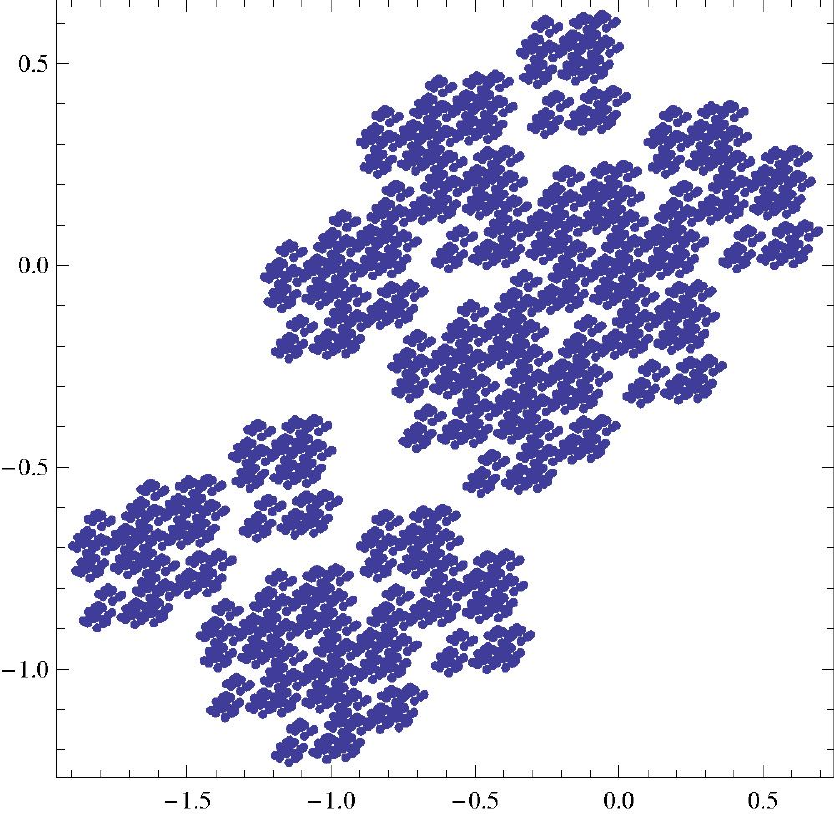}
} \qquad
 \subfigure[$x^2+2x+3$]{
   \includegraphics[width=5cm] {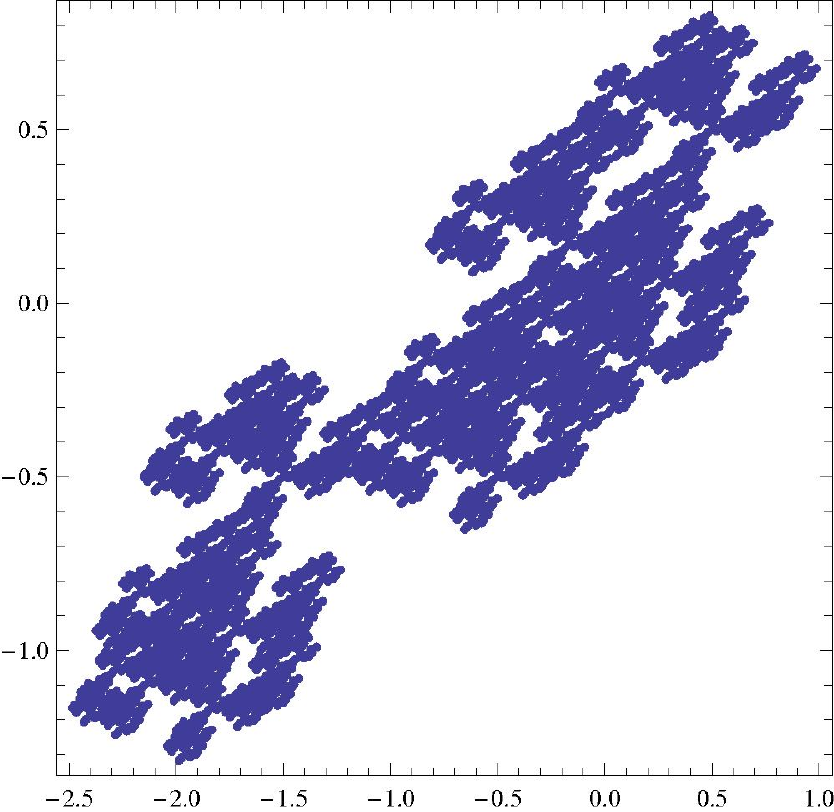}
 }\\
 \subfigure[$x^2+3x+3$]{
 \includegraphics[width=5cm]{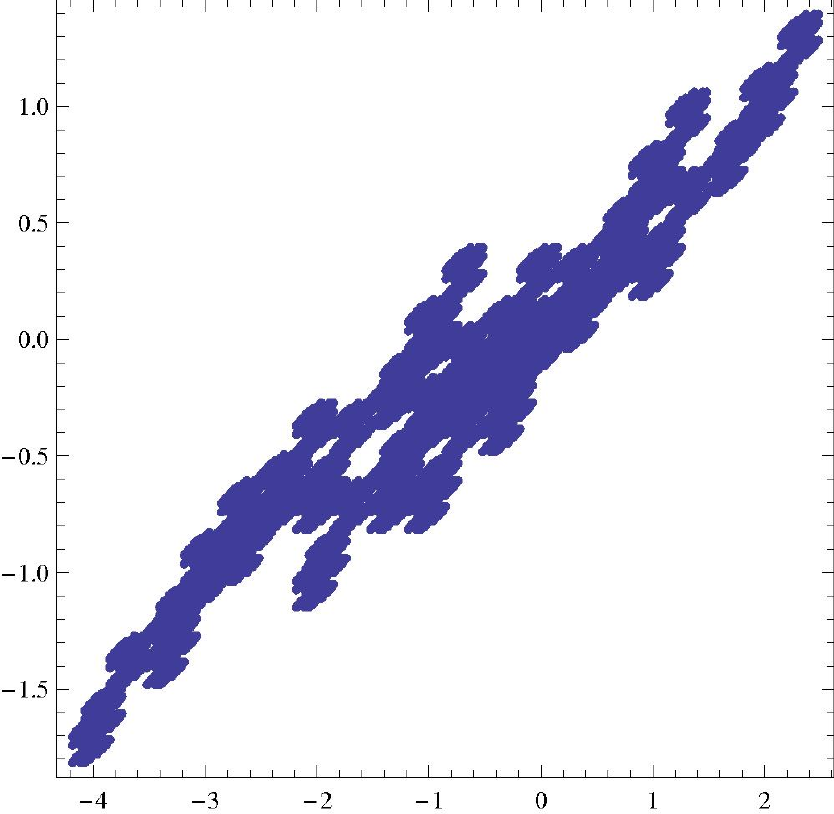}
}
 \qquad
 \subfigure[$x^2+x-3$]{
   \includegraphics[width=5cm]{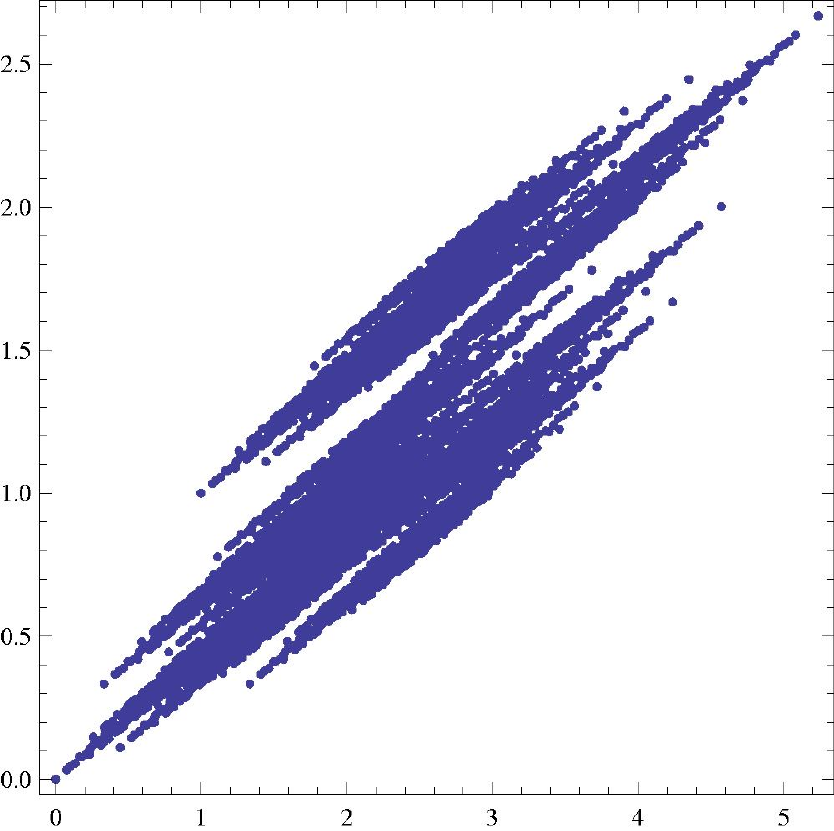}
 }
\caption{When $p\ne 0, m=3$.}\label{figure 2}
\end{figure}

\medskip
\subsection{When $p\ne 0$ and $m\leq 3$} (see  Figure \ref{figure 2} where
we take the vector $v=\left[\begin{array}{rr}
1\\0\end{array}\right]$)

\begin{theorem} Let ${\mathcal{D}}=\{0,1,m\}v$ be a digit set with $2\leq m\in {\mathbb Z}$ and $v\in {\mathbb{R}}^2$ such that
$\{v,Av\}$ is linearly independent. If the characteristic polynomial
of the matrix $A$ is $f(x)=x^2\pm x - 3$ or $x^2\pm 2x+3$ or $
x^2\pm 3x+3$.  Then $T(A,{\mathcal D})$ is a connected set if $m=2
~\text{or}~ 3$. Moreover, $T(A,{\mathcal D})$ is a connected tile if
and only if $m=2$.

If $f(x)=x^2\pm x + 3$. Then $T(A,{\mathcal D})$ is a connected set
if and only if $m=2$.
\end{theorem}

\medskip

\begin{proof}
Since $|\det (A)| =3$ is a prime, $\{v,Av\}$ is independent.
Applying Theorem 3.1 in \cite{KiLa}, we know $T(A,{\mathcal D})$ is
a tile if and only if $\{0,1,m\}$ is a complete set of coset
representatives of ${\mathbb{Z}}_3$, that is, $m=2+3k$ for $k\in
{\mathbb{Z}}$. Moreover, the connectedness of $T$ implies $m\leq 4$
by Theorem \ref{main thm. of special digit sets}. Therefore $m=2$.

\medskip

To show the other parts, by Theorem \ref{main thm. of special digit
sets}, it suffices to verify the cases for $m=3$, where the digit
set becomes ${\mathcal{D}}=\{0,1,3\}v$, and
$\Delta{\mathcal{D}}=\Delta Dv=\{0,\pm 1,\pm 2,\pm 3\}v$.

\medskip

For $f(x)=x^2- x -3$, we have $f(A)=A^2-A-3I=0$ which implies
$v=A^{-1}v+3A^{-2}v \in T-T$, that is,  $T\cap(T+v)\ne \emptyset$.

On the other hand, $0=f(A)(2A-I)=2A^3-3A^2-5A+3I=(2A-3I)(A^2-I)-3A$,
it follows that
$2A-3I=3A^{-1}(I-A^{-2})^{-1}=3\sum_{i=0}^{\infty}A^{-2i-1}$.
Consequently, $$2v=3A^{-1}v+3\sum_{i=1}^{\infty}A^{-2i}v \in T-T.$$
That is, in radix expansion,  $2=0.3\overline{30}$. Hence
$(T+v)\cap(T+3v)\ne \emptyset$. Therefore  $T$ is connected by
Proposition \ref{e-connected prop}.

\medskip

For $f(x)=x^2+2x+3$, we have $f(A)=A^2+2A+3I=0$ which implies
$v=-2A^{-1}v-3A^{-2}v\in T-T$, that is $T\cap(T+v)\ne \emptyset$.

On the other hand, from $0=(A-I)(2A^2+I)(A^2+2A+3I)$, we obtain
$2A^2+2A+3I=(A^3-I)^{-1}(3A^2-3A)$. It follows that
$223=0.\overline{3(-3)0}$ and the radix expansion of $2v$ is
$$2=0.(-2)(-3)\overline{3(-3)0}.$$  Hence
$(T+v)\cap(T+3v)\ne\emptyset$. Therefore $T$ is connected.

\medskip

For $f(x)=x^2+3x+3$, we have $f(A)=A^2+3A+3I=0$. Then
$I=-3A^{-1}-3A^{-2}$, so $v=-3A^{-1}v-3A^{-2}v\in T-T$, that is,
$T\cap(T+v)\ne \emptyset$.

On the other hand, from $0=(2A-I)(A^2+3A+3I)=2A^3+5A^2+3A-3I$, we
obtain $2A^2+3A=3(A+I)^{-1}=3\sum_{i=1}^{\infty}(-1)^{i-1}A^{-i}$.
It follows that  $230=0.\overline{3(-3)}$ and the radix expansion of
$2v$ is
$$2=0.(-3)0\overline{3(-3)}.$$
Hence $(T+v)\cap(T+3v)\ne\emptyset$. Therefore $T$ is connected.

\medskip

For $f(x)=x^2+ x+3$, we can use the same method in the proof of
Theorem \ref{main thm. of special digit sets} to show that
$T(A,{\mathcal{D}})$ is disconnected. Indeed, let $l=\gamma v
+\delta Av=\sum_{i=1}^{\infty}b_iA^{-i}v$, where $b_i\in \Delta
D=\{0,\pm 1,\pm 2,\pm 3\}$, then $T+l$ is a neighbor of $T$. By
(\ref{estimate of gamma and delta 1}), we have estimates
$$|\gamma|<2.64, \quad |\delta|< 1.89.$$

Suppose $(T+v)\cap (T+3v)\ne\emptyset$,
$2v=\sum_{i=1}^{\infty}b_iA^{-i}v$. By (\ref{inductive formula of
gamma_n and delta_n}), it follows that $l_1=\gamma_1 v+\delta_1
Av=-(6+b_2)v-(2+b_1)Av$, where $|6+b_2|<2.64$, which is impossible
since $|b_2|\leq 3$.

Suppose $T\cap (T+3v)\ne \emptyset$, then
$3v=\sum_{i=1}^{\infty}b_iA^{-i}v$. By (\ref{inductive formula of
gamma_n and delta_n}), it follows that $l_1=-(9+b_2)v-(3+b_1)Av$,
where $|9+b_2|< 2.64$, which is also impossible as $|b_2|\leq 3$.
Therefore $T$ is disconnected.

\end{proof}

\end{section}

\bigskip

\begin{section} {\bf General collinear digit set: ${\mathcal{D}}=\{0,1,b\}v, 1<b\in{\mathbb R}.$}

In this section, we discuss the case when $|\det (A)|=3$ and
${\mathcal{D}}=\{0,1,b\}v$, where $1<b\in{\mathbb R}$. By modifying
the proof of Theorem \ref{main thm. of special digit sets}, we prove
Theorem \ref{main results} which is restated in the following.

\medskip

\begin{theorem}
Assume that $A\in M_2({\mathbb{Z}})$ with $|\det (A)| =3$ is
expanding, a digit set ${\mathcal{D}}=\{0,1,b\}v$, where $b>1$ and $
v\in {\mathbb{R}}^2$ such that $\{v,Av\}$ is linearly independent.
Then we have:

 Case 1: $f(x)= x^2\pm x + 3\quad T ~\text{is disconnected \quad
if \quad} b \geq 67/25 ~\text{or}~  b \leq 67/42;$

 Case 2: $f(x)= x^2\pm 2x + 3\quad T ~\text{is disconnected \quad
if \quad}  b \geq 37/10 ~\text{or}~  b \leq 37/27;$

 Case 3: $f(x)= x^2\pm 3x + 3\quad T ~\text{is disconnected \quad
if \quad} b \geq 33/10 ~\text{or}~  b \leq 33/23;$

 Case 4: $f(x)= x^2\pm x - 3\quad T ~\text{is disconnected \quad
if \quad} b >19/5 ~\text{or}~ b<19/14.$
\end{theorem}

\medskip

\begin{proof}
In each case, we mainly consider the situation of  $b\geq 2$ and get
one desired inequality. For otherwise, if $1 <b\leq 2$, then
$b/(b-1)\geq 2$, we replace the original digit set with
${\mathcal{D}}'=\{0,1,b/(b-1)\}v$. Hence the other inequality
follows. If $T$ is connected, then
\begin{eqnarray*}
(b-y)v=\sum_{k=1}^{\infty}b_kA^{-k}v\quad\text{ holds for}\quad y=0
\quad\text{or}\quad 1
\end{eqnarray*}
where $b_k\in \Delta D=\{0,\pm 1,\pm(b-1),\pm b\}$. That is, $T+
(b-y)v$ is a neighbor of $T$.  By making use of the same estimates
 of the neighbors as in the previous section, we prove the theorem case by
case.

\medskip

Case 1: $f(x)= x^2+ x +3$. By (\ref{inductive formula of gamma_n and
delta_n}) and (\ref{estimate of gamma and delta 1}), we obtain
$l_1=-(3(b-y)+b_2)v-(b-y+b_1)Av$ and  $|3(b-y)+b_2|<0.88b$. Since
$|3(b-y)+b_2|\geq 3(b-1)-b=2b-3$, it follows that $b<67/25$. On the
other hand, it yields $b>67/42$ from $b/(b-1)< 67/25$.

\medskip

Case 2: $f(x)= x^2+ 2x +3$. By (\ref{inductive formula of gamma_n
and delta_n}) and (\ref{estimate of gamma and delta 2}), we have
$l_1=-(3(b-y)+b_2)v-(2(b-y)+b_1)Av$, and $|3(b-y)+b_2|< 1.17b$.
Since $|3(b-y)+b_2|\geq 2b-3$, it follows $b<3/0.83<37/10$. We
obtain  $b>37/27$ from $b/(b-1)< 37/10$.

\medskip

Case 3: $f(x)= x^2+ 3x +3$. By using (\ref{inductive formula of
gamma_n and delta_n}) and (\ref{estimate of gamma and delta 3})
repeatedly, we have $l_1=-(3(b-y)+b_2)v-(3(b-y)+b_1)Av$ and
$l_2=(9(b-y)+3b_1-b_3)v+(6(b-y)+3b_1-b_2)Av$,  where
$|9(b-y)+3b_1-b_3|<2.24b$. Since $|9(b-y)+3b_1-b_3|\geq 5b-9$, it
follows $b<33/10$. From $b/(b-1)< 33/10$ we get  $b>33/23$.

\medskip

Case 4: $f(x)= x^2+ x -3$. By  using (\ref{inductive formula of
gamma_n and delta_n}) and (\ref{estimate of gamma and delta 4})
repeatedly, we can get another neighbor
$T+\big(-21(b-y)-12b_1+3b_2-3b_3-b_5\big)v+\big(19(b-y)+7b_1-4b_2+b_3-b_4\big)Av$
and $|19(b-y)+7b_1-4b_2+b_3-b_4|\leq b$. Since
$|19(b-y)+7b_1-4b_2+b_3-b_4|\geq 6b-19$,  it follows that $b\leq
19/5$.  It follows from $b/(b-1)\leq 19/5$ that $b\geq 19/14$.
\end{proof}

\end{section}

\bigskip

\begin{section} {\bf Other results on $|\det (A)|>3$}
In previous sections, we investigated the connectedness of certain
planar self-affine sets $T(A, {\mathcal D})$  generated by an
expanding integer matrix $A$ with $|\det (A)|=3$ and  a
non-consecutive collinear digit set ${\mathcal D}$. Generally, it is
very difficult to determine the connectedness of $T(A, {\mathcal
D})$ for $|\det (A)|>3$ since the structure of $\Delta {\mathcal D}$
can become more complicated. In the section, we try to give some
partial answers.

\medskip

Let $f(x)=x^2+px\pm q$ be the characteristic polynomial of an
expanding integer matrix $A$ where $q\geq 2$. Let $v\in {\mathbb
R}^2$ such that $\{v,Av\}$ is linearly independent and ${\mathcal
D}=Dv$ be a $q$-digit set such that $D=\{0=d_1,d_2,\dots,
d_q\}\subset {\mathbb Z}$ in the increasing order with
$d_{i+1}-d_i=1~\text{or}~2 $ for all $i$, $d_{j+1}-d_j=1~\text{for
at least one}~ j $ and $d_{k+1}-d_k=2~\text{for at least one}~ k $.
We have the following sufficient conditions for this $T(A, {\mathcal
D})$ to be connected. Since we have solved the case when $p=0$
(Theorem \ref{thm of p=0}) and it is known that $T(A, {\mathcal D})$
is connected if and only if $T(-A, {\mathcal D})$ is connected
(Corollary \ref{symmetry of char. poly.}), we may assume $p>0$ in
this section.

\medskip

\begin{theorem} \label{Thm1}
Let $f(x)=x^2+px+ q$ with $q\geq 2$ and $2p>q+2$. Assume $\{0,\pm
1,\pm 2,\dots, \pm (q-1)\}\subset \Delta D$. Then $T$ is connected
if $2p-2\in \Delta D$ and $2q-p\in \Delta D$.

\end{theorem}

\medskip

\begin{proof}
By using $0=f(A)(A-I)=A^3+(p-1)A^2+(q-p)A-qI$, we have

\begin{equation}\label{equ.5_1}
1=0.\overline{(1-p)(p-q)(q-1)}.
\end{equation}
Hence $$T\cap (T+v)\ne \emptyset.$$

Similarly, from $A+(p-1)I=-(q-p+1)(A+I)^{-1}$ we obtain

\begin{equation}\label{equ.5_2}
1=0.(1-p)\overline{[-(q-p+1)](q-p+1)}.
\end{equation}

Adding (\ref{equ.5_1}) and (\ref{equ.5_2}), we get

$$2=0.(2-2p)\overline{(2p-2q-1)(2q-p)(-q)1(p-2)(q-2p+2)}.$$

We can see easily that $|2p-2q-1|,q, 1, |p-2|,|q-2p+2|\leq q-1$. So
if, in addition, $2p-2\in \Delta D$ and $2q-p\in \Delta D$, then it
follows that $$T\cap (T+2v)\ne \emptyset.$$

Therefore $T$ is connected.

\end{proof}

\medskip

\begin{Exa}\textrm{
Let  $f(x)=x^2+5x+6$ and $D=\{0,1,2,4,6,8\}$. Assume ${\mathcal
D}=Dv$ is a digit set such that $\{v,Av\}$ is linearly independent.
We can deduce from Theorem \ref{Thm1} that $T(A,{\mathcal D})$ is
connected. (see  Figure \ref{figure 3}(a) where we take the vector
$v=\left[\begin{array}{rr} 0\\1 \end{array}\right]$)}
\end{Exa}

\medskip

\begin{figure}
 \centering
  \subfigure[$x^2+5x+6$]{
 \includegraphics[width=5cm]{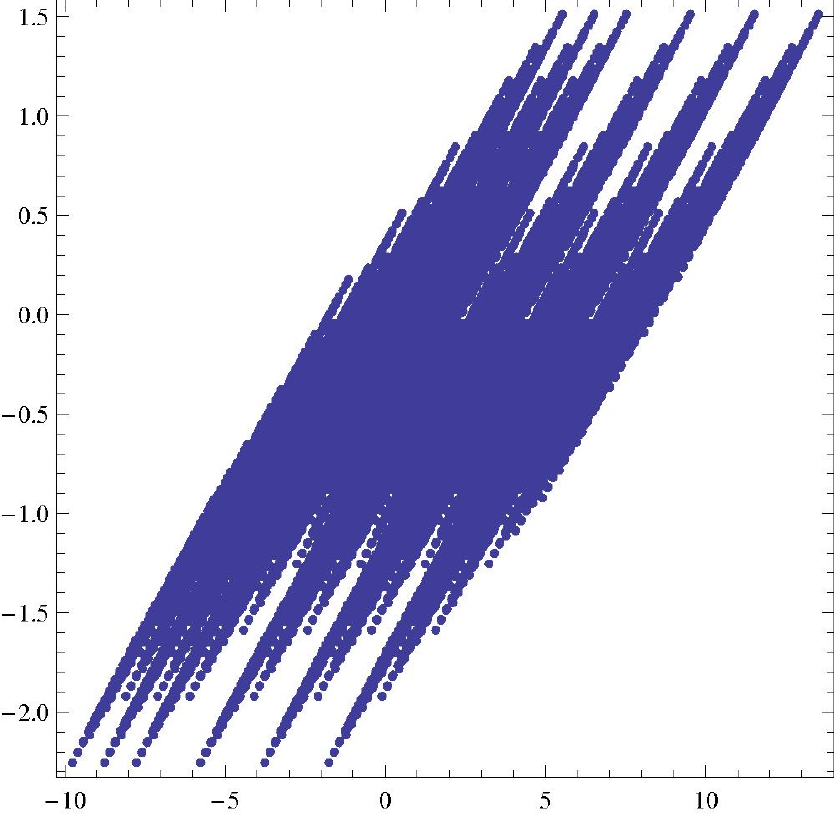}
}
 \qquad
 \subfigure[$x^2+4x-6$]{
   \includegraphics[width=5cm] {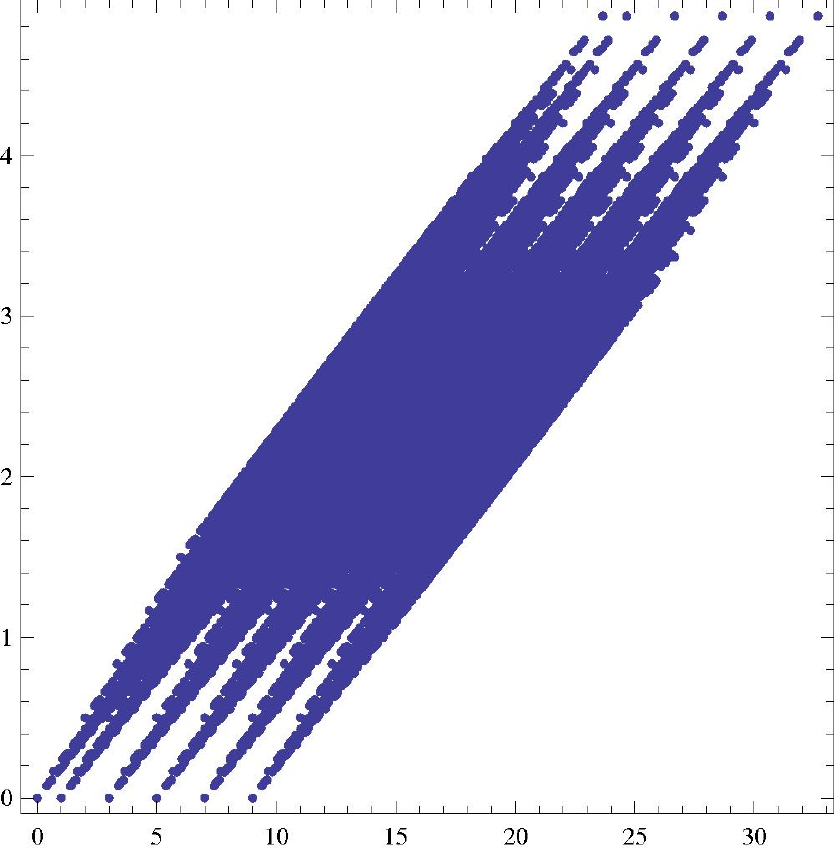}
 }
 \caption{}
 \label{figure 3}
\end{figure}

\begin{theorem} \label{Thm2}
Let $f(x)=x^2+px- q$ with $q\geq 2$ and $2p>q-2$. Assume $\{0,\pm
1,\pm 2,\dots, \pm (q-1)\}\subset \Delta D$. Then $T$ is connected
if $2p+1\in \Delta D$ and $2q-p-2\in \Delta D$.

\end{theorem}

\medskip

\begin{proof}
By using $0=f(A)$, we have

\begin{equation}\label{equ.5_3}
1=0.\overline{(-p)(q-1)}.
\end{equation}
Hence $$T\cap (T+v)\ne \emptyset.$$

Similarly, from $A+(p-1)I=-(q-p-1)(A-I)^{-1}$ we obtain

\begin{equation}\label{equ.5_4}
1=0.[-(p+1)]\overline{(q-p-1)}.
\end{equation}

Adding (\ref{equ.5_3}) and (\ref{equ.5_4}), we get

$$2=0.(-2p-1)\overline{(2q-p-2)(q-2p-1)}.$$

We can see easily that $|q-2p-1|\leq q-1$. So if, in addition,
$2p+1\in \Delta D$ and $2q-p-2\in \Delta D$, then it follows that
$$T\cap (T+2v)\ne \emptyset.$$

Therefore $T$ is connected.

\end{proof}

\medskip

\begin{Exa}\textrm{
Let  $f(x)=x^2+4x-6$ and $D=\{0,1,3,5,7,9\}$. Assume ${\mathcal
D}=Dv$ is a digit set such that $\{v,Av\}$ is linearly independent.
We can deduce from Theorem \ref{Thm2} that $T(A,{\mathcal D})$ is
connected. (see Figure \ref{figure 3}(b) where we take the vector
$v=\left[\begin{array}{rr} 0\\1 \end{array}\right]$) }
\end{Exa}

\end{section}

\bigskip
\bigskip
\noindent {\it Acknowledgements}:  The authors would like to thank
Prof K.S. Lau for suggesting the question and  Prof B. Tan for the
draft of Proposition \ref{Tan's example}. They also thank the
referee for some helpful comments which improve the presentation of
the paper.

\end{document}